\documentclass[a4paper,twoside,11pt]{article}
\usepackage[utf8]{inputenc}
\usepackage[T2A]{fontenc}
\usepackage[english,russian]{babel}
\usepackage{fancyhdr}
\usepackage{newprog1e}
\usepackage{amsfonts,amsmath,amssymb,amsthm}
\usepackage{mathtools}
\usepackage{graphicx}
\usepackage{caption}

\journalnumber{4}
\curyear{2023}

\authorlist{Н.\,A.~Ильтяков и др.}

\titlehead{Об ускоренных покомпонентных методах поиска равновесий в двухстадийной модели}
\headerdef

\udk{519.853.62}
\rubrika{АНАЛИЗ ДАННЫХ}
\dateinput{16.07.2023}

\rusabstr{
Поиск равновесия в двухстадийной модели транспортных потоков сводится к решению специальной негладкой задачи выпуклой оптимизации с двумя группами разных переменных. Для численного решения данной задачи в статье предложено использовать ускоренный блочно-покомпонентный метод Нестерова--Стиха со специальным выбором вероятностей блоков на каждой итерации (одного из двух). Теоретические оценки сложности такого подхода могут заметно улучшать оценки используемых ранее подходов. Однако в общем случае не гарантируют более быстрой сходимости. 
В статье проведены численные эксперименты с предложенными алгоритмами.
\\
Ключевые слова: транспортное моделирование, выпуклая оптимизация, ускоренный покомпонентный градиентный метод
}

\author{
{\bfseries Н.\,A.~Ильтяков\,$^{*}$, М.\,А.~Обозов\,$^{\dagger}$, И.\,М.~Дышлевский\,$^{\ddagger}$, Д.\,В.~Ярмошик\,$^{\it{a},\,\it{b},\,\mathsection}$,
М.\,Б.~Кубентаева\,$^{\it{a},\,\mathparagraph}$,
\\ {\bfseries А.\,В.~Гасников\,$^{\it{a,b,c},\,\|}$,
Е.\,В.~Гасникова \,$^{\it{a},\,**}$}
\\ {\itshape $^{\it{a}}$\,Московский физико-технический институт,}
\\ {\slshape 141701,} {\itshape г. Долгопрудный, Институтский пер.,~9}
\\ {\itshape $^{\it{b}}$\,Институт проблем передачи информации РАН,}
\\ {\slshape 127051,} {\itshape г. Москва, Большой Каретный пер.,~19, стр.~1}
\\ {\itshape $^{\it{c}}$\,Кавказский математический центр Адыгейского гос. университета,}
\\ {\slshape 385000,} {\itshape г. Майкоп, Первомайская ул.,~208}
\\ {\itshape $*$\,e-mail: iltyakov.nik@gmail.com, ORCID: 0009-0004-6090-6469}
\\ {\itshape $\dagger$\,e-mail: obozovmark9@gmail.com, ORCID: 0009-0006-6195-1848}
\\ {\itshape $\ddagger$\,e-mail: igordyslevski@gmail.com}
\\ {\itshape $\mathsection$\,e-mail: yarmoshik.dv@phystech.edu, ORCID: 0000-0003-1912-1040}
\\ {\itshape $\mathparagraph$\,e-mail: kubentay@gmail.com}
\\ {\itshape $\|$\,e-mail: gasnikov@yandex.ru, ORCID: 0000-0002-7386-039X}
\\ {\itshape $**$\,e-mail: egasnikov@yandex.ru}
}}
\title{Об ускоренных покомпонентных методах поиска равновесий в двухстадийной модели равновесного распределения транспортных потоков}

\date{}

\usepackage{pgfplots}
    \usetikzlibrary{
        pgfplots.colorbrewer,
    }
    \pgfplotsset{
        cycle list/.define={my marks}{
            every mark/.append style={solid,fill=\pgfkeysvalueof{/pgfplots/mark list fill}},mark=*\\
            every mark/.append style={solid,fill=\pgfkeysvalueof{/pgfplots/mark list fill}},mark=square*\\
            every mark/.append style={solid,fill=\pgfkeysvalueof{/pgfplots/mark list fill}},mark=triangle*\\
            every mark/.append style={solid,fill=\pgfkeysvalueof{/pgfplots/mark list fill}},mark=diamond*\\
        },
    }\newcommand{\numberthis}{\addtocounter{equation}{1}\tag{\theequation}}

\usepackage{indentfirst} 
\usepackage{pifont}
\usepackage{pgfplots}

\usepackage{graphicx}
\usepackage[colorinlistoftodos]{todonotes}
\usepackage[colorlinks=true, allcolors=blue]{hyperref}

\usepackage{algorithm}
\usepackage{algorithmic}
\floatname{algorithm}{Алгоритм}
\usepackage{varwidth}

\usepackage{wasysym}

\newtheorem{theorem}{Теорема}

\newcommand{\argmin}{\mathop{\arg\!\min}}

\newcommand*{\BREAK}{\textbf{break}}

\def\eps{\varepsilon}

\def \R {\mathbb R}

\newcommand{\demyan}[1]{\textcolor{purple}{#1}}

\DeclareMathOperator{\Mirr}{Mirr}
\DeclareMathOperator{\Grad}{Grad}
\let\Proj\relax
\DeclareMathOperator{\Proj}{Proj}
\DeclareMathOperator{\dom}{dom}

\begin{document}
\maketitle

\section{Введение} \label{section_1}
Многостадийная модель транспортных потоков является основой любого современного пакета транспортного моделирования крупных мегаполисов \cite{Ortuzar2002,Boyles2020}. В основе таких (многостадийных) моделей лежат  задачи выпуклой оптимизации (блоки), последовательное решение (прогонка) которых (по циклу) приближает к искомому равновесному распределению \cite{Ortuzar2002,Boyles2020,gasnikov2020posobie}. Альтернативный путь -- попробовать найти такую общую задачу выпуклой оптимизации, решение которой давало бы искомое равновесие \cite{gasnikov2020posobie}. Альтернативный путь, по-видимому, впервые был предложен в 1976 году С.П. Эванс \cite{Evans1976}. А в современном варианте обоснован в работах А.В. Гасникова и Ю.Е. Нестерова с соавторами \cite{gasnikov2014matmod,Gasnikova2023}. 

Современные численные методы решения задачи выпуклой оптимизации, возникающей при альтернативном подходе базируются на сочетаниях ускоренного универсального метода Нестерова (по группе переменных, отвечающих затратам на ребрах / дорогах) и метода балансировки Брэгмана--Шелейховского по матрице корреспонденций. Этот подход вполне успешно работает на практике \cite{kotlyarova2021,Kubentaeva2023}. Однако, метод совершенно контринтуитивен (<<не физичен>>) по своей сути. В реальной жизни в медленном времени меняется матрица корреспонденций, а в быстром времени под эти изменения подстраиваются затраты на ребрах \cite{Gasnikova2023}. В численном методе все происходит ровно наоборот.

В настоящей статье предложен <<физичный>> численный метод решения той же самой задачи, в котором блок балансировки Брэгмана--Шелейховского заменяется блоком ускоренного метода Нестерова. Точнее говоря, для решения задачи предлагается использовать специальную версию ускоренного блочно-покомпонентного метода Нестерова--Стиха \cite{Nesterov2017} со специальным способом выбора вероятностей блоков (одного из двух).\footnote{Отметим, что хотя статья Нестерова--Стиха вышла в 2017 году, впервые этот методы был доложен ими на конференции, посвященной 80-и летию Б.Т. Поляка в мае 2015 года. Это важно в контексте приоритета Нестерова--Стиха, поскольку почти в то же время, но чуть позже появились близкие работы \cite{gasnikov2016coordinate,Allen-Zhu2016}.}

Отметим, что исходная задача выпуклой оптимизации не является гладкой, поэтому, предварительно потребуется либо ее сгладить, т.е. рассматривать, так называемые, стохастические равновесия в блоке равновесного распределения потоков по путям \cite{baimurzina2019jvm}, либо использовать универсальный вариант метода Нестерова--Стиха \cite{gasnikov2016coordinate}, который сам настраивается на гладкость задачи.

\section{Двухстадийная модель и $\min \min$ задача выпуклой оптимизации}
Данный раздел не является оригинальным во многом базируется на существующей литературе, в частности, на \cite{kotlyarova2021,Gasnikova2023_}.

\textbf{Основные определения и обозначения.}
Будем рассматривать замкнутую транспортную систему, описываемую графом $G = \langle V, E\rangle$, где $V$ -- множество вершин ($|V| = n$), а $E$ -- множество ребер ($|E| = m$). Обозначим ребра графа через $e\in E$. 
Транспортный граф $G$ считается известным.

Часть вершин $O\subseteq V$ (\textit{origin}) являются источниками корреспонденций, а часть стоками корреспонденций $D\subseteq V$ (\textit{destination}). Точнее говоря, вводится множество пар (источник, сток) корреспонденций $OD \subseteq V\bigotimes V$. Сами корреспонденции будем обозначать через $d_{ij}$, где $(i,j)\in OD$. Как правило $|OD|\ll n^2$ \cite{gasnikov2014matmod}. Не ограничивая общности, будем далее считать, что $\sum_{(i,j)\in OD} d_{ij} = 1$. Множество пар $OD$ считается известным. Корреспонденции -- не известны! Однако известны (заданы) характеристики источников и стоков корреспонденций. То есть известны величины $\{l_i\}_{i\in O}$, $\{w_j\}_{j\in D}$
\begin{equation}\label{corr}
    \sum_{j: (i,j)\in OD} d_{ij} = l_i, \quad \sum_{i: (i,j)\in OD} d_{ij} = w_j.
\end{equation}
Заметим, что $\sum_{i\in O} l_i = \sum_{j\in D} w_j = 1$. Условие \eqref{corr} будем также для краткости записывать в виде $d\in (l,w)$.

Обозначим через $\tau_e(f_e)$ -- функцию затрат (например, временных) на проезд по ребру (участку дороги) $e$, если поток автомобилей на этом участке $f_e$. Функции $\tau_e(f_e)$ считаются заданными, например, таким образом: \cite{gasnikov2013book,Patriksson2015,gasnikov2020posobie}
\begin{equation}\label{BPR}
    \tau_e(f_e) = \bar{t}_e\left(1 +\kappa\left(\frac{f_e}{\bar{f}_e}\right)^{\frac{1}{\mu}}\right), 
\end{equation}
где $\bar{t}_e$ -- время прохождения ребра $e$, когда участок свободный (определяется разрешенной скоростью на данном участке), а $\bar{f}_e$ -- пропускная способность ребра $e$ (определяется полосностью: [пропускная способность] $\le$ [число полос] * [2000 авт/час] и характерстиками перекрестков). Считается, что эти характеристики известны \cite{Stabler2020}. Параметр $\mu = 0.25$ -- BPR-функции \cite{Patriksson2015}, но допускается и $\mu\to 0+$ -- модель стабильной динамики \cite{Nesterov2003,gasnikov2013book,gasnikov2014matmod,gasnikov2020posobie,kotlyarova2022}. Параметр $\kappa >0$ также считается заданным.

Полезно также ввести $t_e$ -- (временные) затраты на прохождения ребра $e$. Согласно вышенаписанному $t_e = \tau_e(f_e)$. По этим затратам $t = \{t_e\}_{e\in E}$ можно определить затраты на перемещение из источника $i$ в сток $j$ по кратчайшему пути: 
\begin{equation}\label{T}
T_{ij}(t) = \min_{p \in P_{ij}} T_p(t):= \sum_{e\in E} \delta_{ep}t_e,
\end{equation}
где $p$ -- путь (без самопересечений -- циклов) на графе (набор ребер), $P_{ij}$ -- множество всевозможных путей на графе, стартующих из источника $i$ и заканчивающихся в стоке $j$, $\delta_{ep} = 1$, если ребро $e$ принадлежит пути $p$ и $\delta_{ep} = 0$ -- иначе.

Далее также будет полезен вектор $x = \{x_p\}_{p \in P}$ -- вектор распределения потоков по путям, где $P = \bigcup_{(i,j)\in OD} P_{ij}$. Заметим, что $f_e = \sum_{p} \delta_{ep} x_p$ или в матричном виде $f = \Theta x$, где $\Theta = \|\delta_{ep}\|_{e\in E, p\in P}$.

\textbf{Энтропийная модель расчёта матрицы корреспонденций.}

Под энтропийной моделью расчета матрицы корреспонденций $d(T)$ понимается описанный далее (см. задачу \eqref{Wilson}) способ вычисления набора корреспонденций $\{d_{ij}\}_{(i,j)\in OD}$ по известной матрице затрат $\{T_{ij}(t)\}_{(i,j)\in OD}$. Этот способ заключается в решении задачи энтропийно-линейного программирования (ЭЛП), которую можно понимать, как энтропийно-регуляризованную транспортную задачу
\begin{equation}\label{Wilson}
   \min_{d\in (l,w);d\ge 0} \sum_{(i,j)\in OD} d_{ij}T_{ij}(t) +  \gamma \sum_{(i,j)\in OD} d_{ij}\ln d_{ij},
\end{equation}
где параметр $\gamma >0$ считается известным \cite{Wilson1978,gasnikov2013book,gasnikov2016,gasnikov2020posobie}. Относительно выбора этого параметра, см. \cite{gasnikov2014matmod,gasnikov2020posobie}.

\paragraph{Модели равновесного распределения транспортных потоков по путям.}
Матрица корреспонденций $\{d_{ij}\}_{(i,j)\in OD}$ порождает (вообще говоря, неоднозначно) некий вектор распределения потоков по путям $x$. Неоднозначность заключается в том, что балансовые ограничения, которые возникают на $x \in X(d)$:
\begin{equation*}\label{Xd}
  x\ge 0:\quad \forall (i,j)\in OD \to  \sum_{p\in P_{ij}} x_p = d_{ij},
\end{equation*}
как правило, не определяют вектор $x$ однозначно. Вектор $x$, в свою очередь, порождает вектор потоков на ребрах, $f = \Theta x$, который, в свою очередь, порождает вектор (временных) затрат на ребрах $t(f) = \{\tau_e(f_e)\}_{e\in E}$. На основе последнего вектора уже можно рассчитать матрицу затрат на кратчайших путях $T(t) = \{T_{ij}(t)\}_{(i,j)\in OD}$. Собственно, модель равновесного распределения потоков это формализация \textit{принципа Нэша--Вардропа} о том, что в равновесии каждый водитель выбирает для себя кратчайший путь \cite{gasnikov2013book,Patriksson2015,gasnikov2020posobie}. Другими словами, если для заданной корреспонденции $(i,j)\in OD$ известно, что  (\textit{условие комплиментарности})
\begin{center}
   $x_{p'} >0$, где $p' \in P_{ij}$, то $T_{ij}(t) = \min_{p \in P_{ij}} \sum_{e\in E} \delta_{ep}t_e = \sum_{e\in E} \delta_{ep'}t_e$. 
\end{center}
Задача поиска равновесия сводится, таким образом, к поиску такого вектора $x \in X(d)$, который бы порождал такие затраты $T := T(t(f(x)))$, что выполняется условие комплиментарности. В написанном выше виде искать равновесный вектор $x \in X(d)$ представляется сложной задачей, сводящейся к решению системы нелинейных уравнений. Однако, в данном случае (рассматривается потенциальная игра загрузки) можно свести поиск равновесия к решению задачи выпуклой оптимизации
\begin{equation}\label{Beckman}
  \min_{(f,x): f=\Theta x; x\in X(d)} \sum_{e\in E} \int_{0}^{f_e} \tau_e(z) dz.
\end{equation}
Решение задачи дает модель вычисления вектора потока на ребрах при заданной матрице корреспонденций $f(d)$  \cite{gasnikov2013book,gasnikov2014matmod,Patriksson2015,gasnikov2020posobie}.  

Отметим, что подобно \eqref{Wilson} можно искать не равновесия Нэша--Вардропа, а стохастические равновесия. Это приводит к дополнительному энтропийному слагаемому в \eqref{Beckman} \cite{gasnikov2013book,gasnikov2014matmod,baimurzina2019jvm,gasnikov2020posobie}
\begin{align*}\label{Beckman_stoch}
   \min_{(f,x): f=\Theta x; x\in X(d)} 
   &\sum_{e\in E} \int_{0}^{f_e} \tau_e(z) dz 
   \\
   + &\tilde{\gamma}\sum_{(i,j) \in OD}\sum_{p\in P_{ij}} x_p \ln \left(x_p/d_{ij}\right).\numberthis
\end{align*}
Хотя данной статье планируется все-таки работать с равновесиями Нэша--Вардропа, регулялиризованная задача \eqref{Beckman_stoch} нам понадобится, для того, чтобы сгладить итоговую задачу выпуклой оптимизации, отвечающую двухстадийной модели.

В дальнейшем нам понадобится  (двойственная) переформулировка задач \eqref{Beckman}, \eqref{Beckman_stoch}.  Введем (выпуклые) функции $\sigma_e(f_e) = \int_{0}^{f_e} \tau_e(z) dz$, и обозначим сопряженные к ним функции через $\sigma_e^*(t_e) = \max_{f_e \ge 0} \left\{f_e t_e - \sigma_e(f_e) \right\}$. Например, для BPR-функции \eqref{BPR} $\sigma_e^*(t_e) = \text{const}(\mu)\cdot f_e\cdot(t_e - \bar{t}_e)^{1+\mu}$, при $t_e \ge \bar{t}_e$ и $\text{const}(\mu) \to 1$ при $\mu \to 0+$ \cite{gasnikov2014matmod}.
Тогда (детали см. в \cite{gasnikov2014matmod,gasnikov2020posobie})

\begin{align*}
&\min_{(f,x): f=\Theta x; x\in X(d)} \sum_{e\in E} \sigma_e(f_e) 
\\&= \min_{(f,x): f=\Theta x; x\in X(d)} \sum_{e\in E} \max_{t_e\in \text{dom}\sigma_e^*}\{f_e t_e - \sigma_e^*(t_e)\} 
\\&= 
\max_{\substack{t_e\in \text{dom }\sigma_e^*, \\
e \in E}} \left\{\min_{(f,x): f=\Theta x; x\in X(d)} \sum_{e\in E} f_e t_e\right\}  - \sum_{e\in E}  \sigma_e^*(t_e) 
\end{align*}   
\begin{equation}\label{dual_}
 =   \max_{t_e\in \text{dom }\sigma_e^*, e \in E} \sum_{(i,j)\in OD} d_{ij} T_{ij}(t) - \sum_{e\in E}  \sigma_e^*(t_e). 
 \end{equation}
 Здесь $\text{dom }\sigma_e^*$ означает область определения функции $\sigma_e^*(t_e)$. 

 Примечательно, что задача \eqref{Beckman} имеет самостоятельный и вполне содержательный вывод \cite{Nesterov2003}. 
 А именно, с одной стороны,  $t_e = \tau_e(f_e)$ или
 $f_e = \tau_e^{-1}(t_e) = 
 \frac{d}{dt}\sigma_{e}^{*}(t_e)$ (уравнение состояния транспортного потока:
 чем больше поток по ребру, тем больше времени требуется на прохождение ребра), с другой стороны принцип Нэша--Вардропа (условие комплиментарности), по сути означает, следующее $f \in \partial \sum_{(i,j)\in OD} d_{ij} T_{ij}(t)$, где $\partial$ означает субдифференциал. Осталось заметить, что принцип Ферма (в субдифференциальной форме) для задачи \eqref{dual_} как раз и соответствует двум выписанным соотношениям.
 
 Аналогично можно построить двойственную задачу и к задаче \eqref{Beckman_stoch} (см. также \eqref{T})
 \begin{equation}\label{dual_stoch}
    \max_{t_e\in \text{dom }\sigma_e^*, e \in E} \sum_{(i,j)\in OD} d_{ij} T_{ij}^{\tilde{\gamma}}(t) - \sum_{e\in E}  \sigma_e^*(t_e), 
 \end{equation}
 где $T_{ij}^{\tilde{\gamma}}(t) = -\tilde{\gamma}\ln\left(\sum_{p\in P_{ij}} \exp\left(\frac{-T_p(t)}{\tilde{\gamma}}\right)\right)$.

 \begin{theorem}\label{prop}
      Функция $T_{ij}^{\tilde{\gamma}}(t)$ имеет константу Липшица градиента в $2$-норме равную $H_{ij}/\tilde{\gamma}$, где $H_{ij}$ -- число ребер в самом длинном пути из $P_{ij}$.

      Кроме того, $0 \le T_{ij}^{\tilde{\gamma}}(t) - T_{ij}(t) \le \tilde{\gamma} \ln |P_{ij}|.$
 \end{theorem}
 \textit{Доказательство.} Первое утверждение следует из
 \cite{baimurzina2019jvm}. Второе из \cite{Nesterov2005}.
 $\ocircle$

 \textbf{Двухстадийная модель.} Стандартный способ поиска равновесий в многостадийных транспортных моделях предполагает последовательную прогонку (отрешивание) двух блоков (двух задач) \eqref{Wilson} и \eqref{Beckman}. Из решения \eqref{Wilson} находим зависимость $d(T)$, а из решения \eqref{Beckman} находим зависимость $T(t(f(x(d))))$. Неподвижная точка такой прогонки и будет искомым равновесием. На практике именно такая процедура обычно и реализуется \cite{Ortuzar2002,gasnikov2020posobie,Boyles2020}. Однако, чтобы такая процедура сходилась на практике часто необходимо достаточно удачно выбрать точку старта. Более надежный численный способ поиска равновесия в двухстадийной модели заключается в том, чтобы посмотреть на задачи  \eqref{Wilson} и \eqref{Beckman} (точнее, лучше использовать двойственное представление \eqref{dual_}), и попытаться объединить эти две задачи оптимизации в одну седловую задачу, учитывая их структуры (здесь для краткости соответствующие композитные члены в функционалах задач \eqref{Wilson} и \eqref{dual_} были обозначены $g$ и $h$, ну а общая часть (зависящая от переменных $d$ и $t$) обозначена через $G$):
$$\min_{d\in(l,w)} G(d,T(t)) + g(d),$$
$$\max_{t\in \text{dom }\sigma^*} G(d,T(t)) - h(t).$$
Совместное решение этих двух задач можно найти из решения седловой (выпукло-вогнутой) задачи
$$\min_{d\in(l,w)}\max_{t\in \text{dom }\sigma^*} G(d,T(t)) + g(d) - h(t),$$
которую, в свою очередь, можно переписать как (теорема фон Неймана--Сиона--Какутани)

$$\max_{t\in \text{dom }\sigma^*}\min_{d\in(l,w)} G(d,T(t)) + g(d) - h(t).$$
И уже для последней седловой задачи можно построить двойственную по части переменных $d$. В результате получается задача (вогнутой) оптимизации с двумя блоками переменных: $t$ и блок двойственных переменных для $d$ (множители Лагранжа к ограничениям $d \in (l,w)$). Действительно, возникающую седловую задачу
\begin{align*}
\min_{d\in (l,w);d\ge 0} &\max_{t\in \text{dom }\sigma^*} 
\sum_{(i,j)\in OD} d_{ij} T_{ij}(t) - \sum_{e\in E} \sigma_e^*(t_e) 
\\&+
\gamma \sum_{(i,j)\in OD} d_{ij}\ln d_{ij}.
\end{align*}
можно переписать в виде
\begin{align*}
  \max_{t\in \text{dom }\sigma^*}\quad &\min_{d\in (l,w);\sum_{(i,j)\in OD}d_{ij}=1;d\ge 0} 
  \sum_{(i,j)\in OD} d_{ij} T_{ij}(t) 
  \\&+
  \gamma \sum_{(i,j)\in OD} d_{ij}\ln d_{ij} - \sum_{e\in E} \sigma_e^*(t_e).
\end{align*}

Вспомогательную задачу минимизации можно представить через двойственную к ней:
\begin{align*}\label{TSD}
 &\max_{\substack{t\in \text{dom }\sigma^*\\ (\lambda, \mu)}}
 -\gamma\ln\left(\sum_{(i,j)\in OD}\exp\left(\frac{-T_{ij}(t) + \lambda_i + \mu_j}{\gamma}\right)\right)  
 \\&+
 \langle l,\lambda \rangle + \langle w,\mu \rangle - \sum_{e\in E} \sigma_e^*(t_e) 
 \\&=
 -\min_{\substack{t\in \text{dom }\sigma^*\\ (\lambda, \mu)}} \gamma\ln\left(\sum_{(i,j)\in OD}\exp\left(\frac{-T_{ij}(t) + \lambda_i + \mu_j}{\gamma}\right)\right) 
 \\&- \langle l,\lambda \rangle - \langle w,\mu \rangle + \sum_{e\in E} \sigma_e^*(t_e).\numberthis 
\end{align*}
Обратим внимание, что добавленное по $d$ ограничение $\sum_{(i,j)\in OD}d_{ij}=1$ тавтологично, поскольку следует из $d\in(l,w)$. Тем не менее, удобнее его добавить, чтобы при взятии $\min$ получалась равномерно гладкая функция (типа softmax), а не сумма экспонент, имеющая неограниченные константы гладкости \cite{gasnikov2020posobie}. Множители $\lambda$ и $\mu$ являются двойственными множителями (множителями Лагранжа) к ограничениям $d\in(l,w)$ (см. \eqref{corr}), которые заносятся в функционал (ограничения $\sum_{(i,j)\in OD}d_{ij}=1;d\ge 0$ не заносятся в функционал). Заметим, что если $(t,\lambda,\mu)$ -- решение задачи \eqref{TSD}, то\footnote{Здесь $C_{\lambda}$ и $C_{\mu}$ -- произвольные числа.} $\left(t,\lambda + (C_{\lambda},...,C_{\lambda})^T,\mu+(C_{\mu},...,C_{\mu})^T\right)$ -- также будет решением задачи, т.е. решение задачи \eqref{TSD} не единственное \cite{gasnikov2020posobie}. Заметим также, что, зная $(\lambda,\mu)$, можно посчитать матрицу корреспонденций \cite{gasnikov2020posobie}:
\begin{equation}\label{d}
    d_{ij}(\lambda,\mu)=\frac{\exp\left(\frac{-T_{ij}(t) + \lambda_i + \mu_j}{\gamma}\right)}{\sum_{(k,l)\in OD}\exp\left(\frac{-T_{kl}(t) + \lambda_k + \mu_l}{\gamma}\right)}.
\end{equation}

Записывая в обратном порядке цепочку равенств \eqref{dual_}, получаем, что прямая задача для двойственной задачи \eqref{TSD} имеет вид
\begin{align*}\label{eq:primal}
   \min_{\substack{(f,x): f=\Theta x; x\in X(d), \\d \in (l, w), d\geq 0}} \bigg\{& P(f, d) =  
   \sum_{e\in E} \sigma_e(f_e) 
   \\&+ 
   \gamma \sum_{(i,j)\in OD} d_{ij}\ln d_{ij}
  \bigg\} \numberthis.
\end{align*}

Для решения (двойственной) задачи выпуклой оптимизации (но, вообще говоря, негладкой, поскольку функции $T_{ij}(t)$ -- негладкие)  можно использовать субградиентные методы. А именно, субградиент (далее обозначаем (супер-)субградиент таким же символом, как и градиент $\nabla$) целевого функционала по $t$ (стоящего под минимумом) \eqref{TSD} можно посчитать по формуле Демьянова--Данскина (см., например, \cite{gasnikov2020posobie}):

\begin{align*}\label{grad}
   &-\sum_{(i,j)\in OD} d_{ij}(\lambda,\mu)\nabla T_{ij}(t) + f 
   \\&=
   -\sum_{(i,j)\in OD} d_{ij}(\lambda,\mu)\nabla T_{ij}(t) +  \{\tau_e^{-1}(t_e)\}_{e\in E},\numberthis
\end{align*}
где, как и ранее, $\tau_e^{-1}$ -- обратная функция к $\tau_e$.

\section{Численные методы решения задачи поиска равновесия в двухстадийной модели}
Перепишем задачу \eqref{TSD} (убрав знак минус) в немного более компактном виде
\begin{align*}\label{TSD_}
 \min_{t\in \text{dom }\sigma^*; (\lambda, \mu)}
\bigg\{& D(t, \lambda, \mu) =
\\&  
\text{softmax}_{\gamma}\left(\left\{-T_{ij}(t) + \lambda_i + \mu_j\right\}_{(i,j)\in OD}\right)
 \\&-
 \langle l,\lambda \rangle - \langle w,\mu \rangle + \sum_{e\in E} \sigma_e^*(t_e) \bigg\}, 
 \numberthis
\end{align*}
где через $D(t, \lambda, \mu)$ мы обозначили двойственную функцию со знаком минус.

Данная задача имеет ярко выраженную блочную $\min \min$ структуру. В частности, задача 
\begin{align}\label{eq:D_t}
\min_{(\lambda, \mu)}  D(t, \lambda, \mu) =: D(t)  
\end{align}
по блоку переменных $(\lambda,\mu)$
может быть довольно эффективно решена методом альтернированных направлений \cite{Peyre2019} (балансировка Брэгмана--Шелейховского--Синхорна) или его ускоренным вариантом \cite{Guminov2019}  с 
достаточно высокой точностью \cite{Tupitsa2022}. 
При этом по $t$ задача негладкая и ее предлагается решать (см., например, \cite{Kubentaeva2023}) универсальным ускоренным методом Нестерова \cite{Nesterov2015,gasnikov2018jvm}, который сам настраивается на гладкость задачи. При таком подходе число итераций определяется, по сути, только гладкостью (Липшицевостью) задачи \eqref{TSD_}. А стоимость итерации определяется согласно \eqref{grad} сложностью метода альтернированных направлений (с учетом возможностей теплого старта метода) и (параллельным!) вычислением (с помощью алгоритма Дейкстры или Беллмана--Форда) субградиента по $t$. Такой подход эффективен, если сложностью метода альтернированных направлений меньше сложности вычисления субградиента по $t$ \cite{gasnikov2020posobie}. Однако такой метод, как уже отмечалось во введении, контринтуитивен, поскольку не соответствует реальному процессу формирования равновесия в двухстадийной модели <<в жизни>> \cite{Gasnikova2023}. 

Альтернативным способом решения задачи \eqref{TSD_} является двойственное сглаживание (по Ю.Е. Нестерову \cite{Nesterov2005}) функций $T_{ij}(t) \to T_{ij}^{\tilde{\gamma}}(t) = -\text{softmax}_{\tilde{\gamma}}\left(\left\{-T_p(t)\right\}_{p\in P_{ij}}\right)$ -- см. теорему \ref{prop}. Итак, рассмотрим задачу
 \begin{align*}\label{TSD_smth}
 &\min_{t\in \text{dom }\sigma^*; (\lambda, \mu)} \text{softmax}_{\gamma}\left(\left\{-T^{\tilde{\gamma}}_{ij}(t) + \lambda_i + \mu_j\right\}_{(i,j)\in OD}\right) 
 \\&- \langle l,\lambda \rangle - \langle w,\mu \rangle + \sum_{e\in E} \sigma_e^*(t_e). \numberthis
\end{align*}
\begin{theorem}\label{smoothing} Решение выпуклой задачи \eqref{TSD_smth} с точностью $\varepsilon/2$ по функции будет решением выпуклой задачи \eqref{TSD} с точностью $\varepsilon$ по функции, если
$$\tilde{\gamma}\le \frac{\varepsilon}{ \max_{(i,j) \in P_{ij}} \ln |P_{ij}|}.$$
Более того, первое слагаемое (softmax) целевой функции в \eqref{TSD_smth} имеет константу Липшица в $2$-норме по $t$ равную $2\sqrt{H}$ ($H = \max_{(i,j)\in OD} H_{ij}$) и константу Липшица градиента в $2$-норме по $t$ равную $H/\min\{\gamma,\tilde{\gamma}\} = H/\tilde{\gamma}$.
При это по $(\lambda,\mu)$ аналогичные константы будут равны, соответственно, $2$ и $1/\gamma$.
 \end{theorem}
\textit{Доказательство}. Первое утверждение следует из теоремы~\ref{prop} и свойства функции softmax -- константа Липишица в $\infty$-норме равна $1$, см. также раздел 2.2.7 \cite{gasnikov2020posobie}. Второе утверждение следует из Замечания 1.4.4 \cite{gasnikov2020posobie} и \cite{Nesterov2005,Guminov2019}.  $\ocircle$

Теорема \ref{smoothing} дает инструмент для использования ускоренного композитного\footnote{Композитом будут все слагаемые в \eqref{TSD_smth} кроме первого нетривиального слагаемого -- softmax. Для ряда конкретных функций $\sigma_e^*(t_e)$ (в том числе, порожденной BPR-функциями затрат) композитность важна, потому что $\sigma_e^*(t_e)$ может иметь неограниченную константу Липшица производной.} блочно-покомонентного метода \cite{gasnikov2016coordinate,Nesterov2017,Allen-Zhu2016} с двумя блоками ($t$ и $(\lambda,\mu)$) или тремя блоками ($t$, $\lambda$ и $\mu$). Строго говоря, методы, приведенные во всех трех работах, не являются композитными. Однако нужное обобщение получается стандартным образом \cite{gasnikov2018jvm}. Отметим, что у всех методов из перечисленных статей есть гиперпараметр, который разумно выбирать таким образом, чтобы вероятности выбора соответствующих блоков были пропорциональны отвечающим этим блокам константам Липшица градиента целевой функции по группе переменных, входящих в блок. Такой способ выбора вероятностей приводит к полному расщеплению оракульных сложностей по блокам \cite{Kovalev2022}. В теоретическом плане это может привести к существенному сокращению объема вычислений по сравнению с контринтуитивным подходом.\footnote{Тут сложно говорить точно, потому что сложность универсального ускоренного метода варьируется от сходимости обычного ускоренного метода в гладком случае (оптимистичный сценарий) до скорости сходимости субградиентного метода в негладком случае (пессимистичный сценарий) \cite{gasnikov2018jvm}. В частности, при пессимистичном сценарии число вычислений $\nabla T_{ij}(t)$ пропорционально $\sim \varepsilon^{-2}$, в то время как в новом подходе аналогичных по сложности вычислений $\nabla T_{ij}^{\tilde{\gamma}}(t)$ потребуется $\sim\max_{(i,j) \in P_{ij}} \ln |P_{ij}|\varepsilon^{-1}$, где $\varepsilon$ -- желаемая по функции точность решения задачи.} 

Недостатком нового подхода является возникновение дополнительного фактора $\max_{(i,j) \in P_{ij}} \ln |P_{ij}|$ в оценке сложности числа вызовов оракула $\nabla T_{ij}^{\tilde{\gamma}}(t)$. Для сетей типа Манхэттенских этот множитель может быть порядка корня из числа ребер графа транспортной сети, что может нивелировать выигрыш от возможно лучшей зависимости оценки сложности от $\varepsilon$. Также вычисление $\nabla T_{ij}^{\tilde{\gamma}}(t)$ хотя и сопоставимо в теории по сложности со сложностью вычисления $\nabla T_{ij}(t)$ \cite{baimurzina2019jvm}, но требует во многом самостоятельной реализации, т.е. тут сложнее воспользоваться готовыми пакетами \cite{gasnikov2020posobie}. Зато, как несложно проверить, новый подход будет соответствовать реальной хронологии событий, то есть будет <<физичным>> (<<интуитивно>> более понятным).

Перечисленные недостатки могут быть устранены, если, следуя \cite{gasnikov2016coordinate}, решать исходную задачу \eqref{TSD_smth} универсальным ускоренным блочно-покомпонентным методом. Строго говоря, теоретического обоснования такого метода на данный момент не известно, однако по аналогии с универсальным ускоренным методом Нестерова можно написать и соответствующую блочно-покомпонентную версию \cite{gasnikov2016coordinate,Nesterov2017,Allen-Zhu2016}. Полученный в результате метод остается <<физичным>>, но наследует хорошее свойство адаптивной настройки на параметры гладкости от <<нефизичного>> метода. 

В заключение этого раздела отметим, что все перечисленные подходы по своей структуре являются прямо-двойственными. Это означает, что существует эффективный способ вычисления $d$ и $f$ по генерируемым последовательностям точек $t, \lambda, \mu$ \cite{gasnikov2020posobie,Kubentaeva2023}.

Подробнее о рассмотренных методах и их работе на практике будет написано в следующем разделе.

\section{Основные алгоритмы и сравнительный анализ их работы}

Для проведения сравнительного анализа были рассмотрены и имплементированы несколько различных алгоритмов оптимизации, а именно: 
 \begin{itemize}
    \item USTM \cite{GasinkovNesterov2016}~-- ускоренный метод подобных треугольников;
    \item USTM + Sinkhorn(балансировка Брэгмана--Шелейховского--Синхорна, описанная выше) \cite{Yarmoshik+Meruza USTM}~-- ускоренный метод подобных треугольников для функции $D(t)$ \eqref{eq:D_t} и отрешивание задачи по $\lambda$ и $\mu$ с помощью алгоритма Синхорна;
    \item ACRCD* \cite{ACRCD}~-- ускоренный, адаптивный покомпонентный метод оптимизации.
\end{itemize}

\noindent
Стоит отметить, что, классический USTM, хотя ранее и не применялся в таком виде для поиска равновесий в транспортной сети, ввиду своей естественности был взят именно в качестве <<бейзлайна>>.

\textbf{\textit{ Алгоритм USTM }}

Листинг алгоритма USTM (Universal Method of Similar Triangles) для минимизации произвольной выпуклой, непрерывной по Липшицу функции $\Phi$ приведён в Алгоритме \ref{alg::univ_triangles}. 
При этом мы использовали следующие обозначения
\begin{equation*}
\phi_0(t) = \frac{1}{2} \|t - t^0\|_2^2,
\end{equation*}
\begin{align*}
&\phi_{k+1}(t) = \phi_k(t) 
\\&+ \alpha_{k+1} \left[{\Phi}(y^{k+1}) 
+ 
\left\langle{\nabla} \Phi(y^{k+1}), t - y^{k+1} \right\rangle \right].
\end{align*}

\begin{algorithm}[h]
    \caption{Универсальный метод подобных треугольников}
    \label{alg::univ_triangles}
    \begin{algorithmic}[1]
    \REQUIRE $L_0 > 0$, стартовая точка $t^0$, точность $\eps > 0$
    \STATE $u^0 \coloneqq t^0$, $A_0 \coloneqq 0$, $k \coloneqq 0$
    \REPEAT
        \STATE $L_{k+1} \coloneqq L_k / 2$ 
        \WHILE{\TRUE}
            \STATE $\alpha_{k+1} \coloneqq \frac{1}{2 L_{k+1} } + \sqrt{\frac{1}{4 L_{k+1}^2} + \frac{A_k}{L_{k+1}} }, \quad A_{k+1} \coloneqq A_k + \alpha_{k+1}$
            \STATE $y^{k+1} \coloneqq \frac{\alpha_{k+1} u^k + A_k t^k}{A_{k+1}}$
            \STATE $u^{k+1} \coloneqq \argmin\limits_{t \in \text{dom} \Phi} \phi_{k+1}(t)$
            \STATE $t^{k+1} \coloneqq \frac{\alpha_{k+1} u^{k+1} + A_k t^k}{A_{k+1}}$
            \IF{${\Phi}(t^{k+1}) \le {\Phi}(y^{k+1}) + \left\langle {\nabla} \Phi(y^{k+1}), t^{k+1} - y^{k+1} \right\rangle 
                + \frac{L_{k+1}}{2} \|t^{k+1} - y^{k+1}\|_2^2 + \frac{\alpha_{k+1}}{2 A_{k+1}} \eps$}
                \STATE \BREAK
            \ELSE
                \STATE $L_{k+1} \coloneqq 2 L_{k+1}$
            \ENDIF
        \ENDWHILE
        \STATE $k \coloneqq k + 1$
    \UNTIL{Критерий остановки не выполнен}
    \end{algorithmic}
\end{algorithm}

В этой статье алгоритмом \textit{USTM} обозначается Алгоритм \ref{alg::univ_triangles}, применённый к функции $D(t, \lambda, \mu)$.
Градиент функции $D(t, \lambda, \mu)$ по переменной $t$ вычисляется по формуле \eqref{grad}, градиент по переменным $\lambda, \mu$ вычисляется напрямую.

\textbf{\textit{
Алгоритм USTM-Sinkhorn
}}

Алгоритмом \textit{USTM-Sinkhorn} мы называем Алгоритм~\ref{alg::univ_triangles}, применённый к функции $D(t)$.
Основное отличие USTM-Sinkhorn от классического USTM~-- использование Алгоритма~\ref{alg:sink} (Синхорна) для решения задачи по $\lambda$ и $\mu$.
Градиент функции $D(t)$ при этом считается так же по формуле \eqref{grad}, куда в качестве $\lambda$ и $\mu$ подставляются значения, найденные "Синхорном".

\begin{algorithm}[h]
    \caption{Алгоритм Синхорна}
    \label{alg:sink}
    \begin{algorithmic}[1]
    \STATE $\lambda^0 \coloneqq \vec 0, \mu^0 \coloneqq \vec 0$
    \STATE $k \coloneqq 0$
    \REPEAT
    	\IF{$k$ mod $2=0$}
            \STATE $\lambda^{k+1} \coloneqq \argmin_\xi D(t, \xi, \mu^k)$
            \STATE $\mu^{k+1} \coloneqq \mu^{k}$
        \ELSE
            \STATE $\lambda^{k+1} \coloneqq \lambda^{k}$
            \STATE $\mu^{k+1} \coloneqq \argmin_\xi D(t, \lambda^k, \xi)$
        \ENDIF
        \STATE $k \coloneqq k + 1$
    \UNTIL{Критерий остановки не выполнен}
\end{algorithmic}
\end{algorithm}
Одна итерация алгоритма Синхорна имеет арифметическую сложность, примерно равную арифметической сложности вычисления градиентов двойственной функции по $\lambda$ и $\mu$.

\textbf{\textit{
Алгоритм ACRCD*
}}

Опишем вариант ускоренного покомпонентного метода (Accelerated by Coupling Randomized Coordinate Descent – ACRCD) на базе специального каплинга покомпонентных
вариантов ПГМ (Grad) и МЗС (Mirr) (${x_0 = y_0 = z_0}$), см. Алгоритм~\ref{alg:acrcd}.

Определим
\begin{align*}
    \Grad_i(x) &= \Proj_{\dom f}\left( x - \frac{1}{L_i} \nabla_i f(x)\right),
    \\
    \Mirr_z(\xi) &= \Proj_{\dom f} \left(\left\{ z_i - \frac{1}{L_i}\xi_i \right\}_{i=1}^n \right),
\end{align*}
где $\Proj_Q(x)$~-- евклидова проекция вектора $x$ на множество $Q$, $\nabla_i$~-- градиент по $i$-му блоку переменных (компонент), $L_i$~-- константа Липшица градиента $\nabla_i f(x)$, $n$~-- количество блоков переменных.
В нашем случае $n=2$, блоками переменных являются $t$ и $(\lambda, \mu)$.

В ACRCD* количество вычислений градиентов по каждому блоку обратно пропорционально гладкости функции по этому блоку переменных.
Так как функция $D(t, \lambda, \mu)$ негладкая по блоку переменных $t$ и гладкая по $(\lambda, \mu)$, при правильной инициализации констант $L_i$ метод будет делать больше вычислений градиентов по $t$, чем по $\lambda, \mu$ (см Раздел~\ref{section_5}), что противоположно поведению алгоритма USTM+Sinkhorn, и в то же время согласуется с механизмом установления равновесия в реальной жизни: потоки на путях $f$ (и времена проезда $t$) меняются чаще, чем корреспонденции $d$ (и ``привлекательности'' районов $\lambda, \mu$).

\begin{algorithm}[h]
    \caption{ACRCD*}
    \label{alg:acrcd}
    \begin{algorithmic}[1]
    \REQUIRE  $\tau \in [0,1]$, $\alpha > 0$, $L_i > 0, i \in [1,..,n]$, стартовая точка $x^0$.
    \STATE $y^0 \coloneqq x^0$, $z^0 \coloneqq x^0$.
    \STATE $S \coloneqq \sum_{i=1}^n L_i^{1/2}$
    \REPEAT
        \STATE $\tau_k \coloneqq \frac{2}{k+2}$
        \STATE $x_{k+1} \coloneqq \tau_k z_k + (1- \tau_k)y_k$
        \STATE Выбрать $i_{k+1} \in [1,.., n]$ с вероятностью $P\left(\left\{i_{k+1} = i\right\}\right) = L_i^{1/2}/S$
        \STATE $y_{k+1} \coloneqq \Grad_{i_{k+1}}(x_{k+1})$
        \STATE $\alpha_{k+1} \coloneqq \frac{k+2}{2S^2}$
        \STATE $z_{k+1} \coloneqq \Mirr_{z_k}\left(\alpha_{k+1} S \nabla_{i_{k+1}}f(x_{k+1})\right)$
    \UNTIL{Критерий остановки не выполнен}
    \end{algorithmic}
\end{algorithm}
\noindent

\textbf{Численные эксперименты}

Для сравнения эффективности на практике, перечисленные алгоритмы применялись к задаче поиска равновесия двухстадийной модели на транспортных сетях городов \textit{Anaheim} и \textit{Sioux-Falls} \cite{Stabler2020}.
Параметр $\gamma$ энтропийной модели (формула \eqref{Wilson} и последующие) был принят равным $10$, что по порядку соответствует типично используемым на практике значеням (это важно, потому что скорость сходимости алгоритма Синхорна сильно зависит от величины этого параметра). 
Например, в таком случае энтропийная модель эквивалентна гравитационной модели расчета матрицы корреспонденций с функцией тяготения $\exp(-0.1T_{ij})$ \cite{Aliev2015}, если затраты на межрайонные перемещения $T_{ij}$ выражены в минутах. 



Для оценки сходимости алгоритмов использовалась следующая метрика:
\newcommand{\lamustar}{{\begin{pmatrix} \lambda^*\\ \mu^* \end{pmatrix}}}
\begin{align*}\label{eq:metric}
\Bigg(2\left\| \begin{pmatrix}\nabla_\lambda D(t, \lambda, \mu) \\ \nabla_\mu D(t, \lambda, \mu)  \end{pmatrix}\right\|  \left\|\lamustar \right\| 
\\
 + \max\left\{\text{gap}(t, \lambda, \mu), 0\right\}\Bigg)^{-1}\numberthis,\end{align*}
где $\lamustar$~-- оптимальные значения $\lambda$ и $\mu$ с минимальной нормой(в экспериментах использовались значения $\lambda$ и $\mu$, полученные на последней итерации),
а $\text{gap}$~-- зазор двойственности~-- определяется как разность значений прямой \eqref{eq:primal} и двойственной \eqref{TSD_} функций:
\begin{equation*}
    \text{gap}(t, \lambda, \mu) = P(f, d) + D(t, \lambda, \mu).
\end{equation*}

Такая метрика позволяет учесть оба признака неоптимальности решения: нарушение ограничений и разницу между значением функции в текущей точке и её оптимальным значением, которая оценивается сверху зазором двойственности.
Введение масштабирующего множителя $2\left\|\lamustar \right\|$ у слагаемого, равного норме невязки по ограничениям, позволяет сравнивать между собой невязки по значению функции и по нарушению ограничений. 
Обоснованием такого выбора множителя является прямодвойственный анализ метода градиентного спуска, см (4.13) в \cite{GasnikovUGD}.
\newline
\noindent

Результаты экспериментов приведены на Рисунках \ref{fig:convergence}, \ref{fig:sioux_falls_convergence}, также доступен код вычислительных экспериментов \cite{Code}.

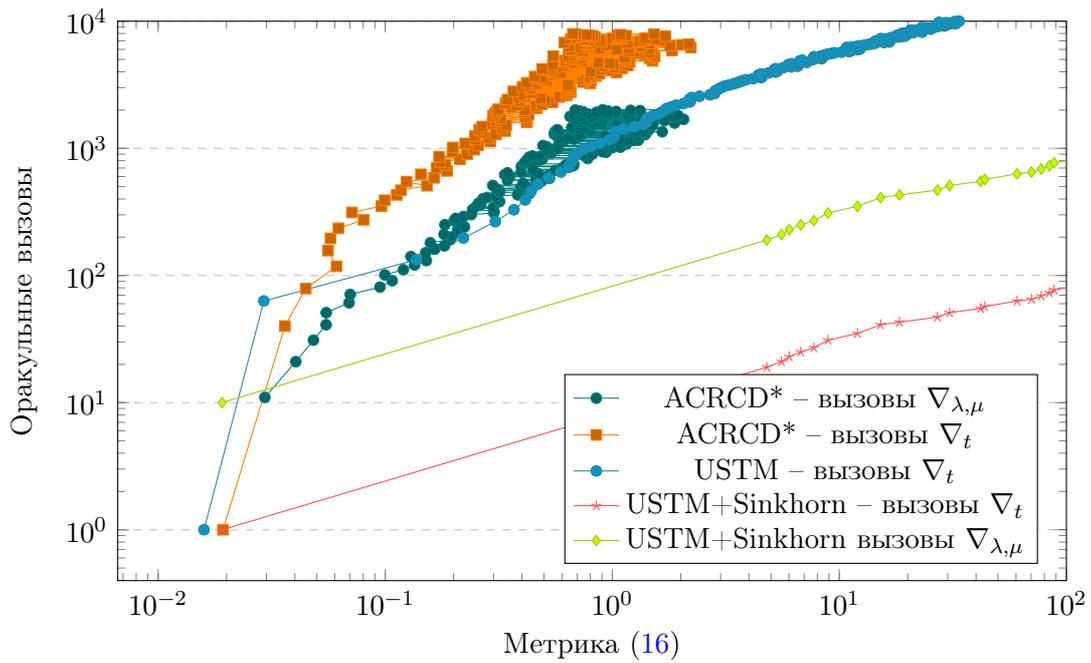
\begin{figure*}
    \centering
        \begin{tikzpicture}
    \begin{axis}
    [
        cycle list name=exotic,
        xlabel={Метрика \eqref{eq:metric}},
        ylabel={Оракульные вызовы},
        xmin=0, xmax=100,
        ymin=0, ymax=10000,
        width=7cm,
        height=9cm,
        xmode=log,
        ymode=log,
        width=400,
        legend pos=south east,
        ymajorgrids=true,
        grid style=dashed,
    ]

    \addplot table [x=x, y=y, col sep=comma] {Acrcd_custom_metric_for_count_la_mu_cals.csv};
    \addplot table [x=x, y=y, col sep=comma] {Acrcd_custom_metric_for_count_t_cals.csv};
    \addplot table [x=x, y=y, col sep=comma] {custom_metric_USTM.csv};
    \addplot table [x=x, y=y, col sep=comma] {custom_metric_USTM_sinkhorn_t_calls.csv};
    \addplot table [x=x, y=y, col sep=comma] {custom_metric_USTM_sinkhorn_sinkhorn_calls.csv};=
    
    \legend{
    {ACRCD* -- вызовы $\nabla_{\lambda, \mu} $},
    {ACRCD* -- вызовы $\nabla_{t}$},
    {USTM -- вызовы $\nabla_{t}$},
    {USTM+Sinkhorn -- вызовы $\nabla_{t}$},
    {USTM+Sinkhorn вызовы $\nabla_{\lambda, \mu}$}
    }
    \end{axis}
    
    \end{tikzpicture}
    \caption{Сравнение алгоритмов(Транспортная задача, граф <<Anaheim>>)}
    \label{fig:convergence}
\end{figure*}
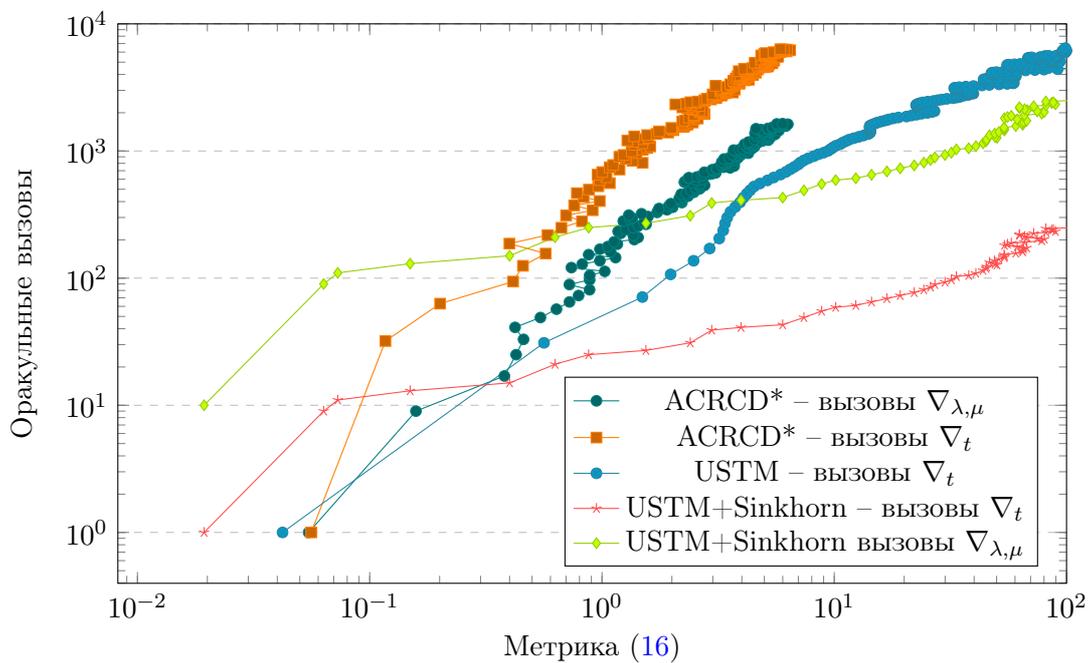
\begin{figure*}
    \centering
        \begin{tikzpicture}
    \begin{axis}
    [
        cycle list name=exotic,
        xlabel={Метрика \eqref{eq:metric}},
        ylabel={Оракульные вызовы},
        xmin=0, xmax=100,
        ymin=0, ymax=10000,
        xmode=log,
        width=7cm,
        height=9cm,
        ymode=log,
        width=400,
        legend pos=south east,
        ymajorgrids=true,
        grid style=dashed,
    ]

    \addplot table [x=x, y=y, col sep=comma] {Acrcd_custom_metric_sufalls_for_count_la_mu_cals.csv};
    \addplot table [x=x, y=y, col sep=comma] {Acrcd_custom_metric_sufalls_for_count_t_cals.csv};
    \addplot table [x=x, y=y, col sep=comma] {custom_metric_sufalls_USTM.csv};
    \addplot table [x=x, y=y, col sep=comma] {custom_metric_sufalls_USTM_sinkhorn_t_calls.csv};
    \addplot table [x=x, y=y, col sep=comma] {custom_metric_sufalls_USTM_sinkhorn_sinkhorn_calls.csv};=

    \legend{
    {ACRCD* -- вызовы $\nabla_{\lambda, \mu} $},
    {ACRCD* -- вызовы $\nabla_{t}$},
    {USTM -- вызовы $\nabla_{t}$},
    {USTM+Sinkhorn -- вызовы $\nabla_{t}$},
    {USTM+Sinkhorn вызовы $\nabla_{\lambda, \mu}$}
    }
    \end{axis}
    
    \end{tikzpicture}
    \caption{Сравнение алгоритмов(Транспортная задача, граф <<Sioux-Falls>>)}
    \label{fig:sioux_falls_convergence}
\end{figure*}

По Рисунку \ref{fig:convergence} видно, что алгоритм USTM+Sinkhorn показал наилучший результат, хотя ему требуются дополнительный итерации, а алгоритм ACRCD*, хоть и обладает нужной <<физичностью>>, но даже не показывает себя лучше <<бейзлана>>. Тест в данном случае проводился на крупном графе <<Anaheim>>, для дополнительной проверки результатов мы провели еще один тест (Рисунок \ref{fig:sioux_falls_convergence}) на маленьком графе <<Sioux-Falls>>. Результаты ожидаемо подтвердились.

Доминирование алгоритма USTM+Sinkhorn мы связываем с тем, что в результате минимизации $D(t, \lambda, \mu)$ по переменным $\lambda, \mu$ получается более гладкая функция $D(t, \lambda, \mu)$ (общее свойство $\min \min$ задач).
Оказалось, что в транспортной задаче положительный эффект от сглаживания с лихвой окупает дополнительные затраты на отрешивание подзадачи по $\lambda, \mu$ алгоритмом Синхорна, хотя этот результат нельзя было предсказать без проведения вычислительных экспериментов.
 
\section{Приложение. О расщеплении оракульных сложностей для негладко-гладких выпуклых задач со структурой $\min \min$}\label{section_5}
Побочным продуктом исследований, проведенных в данной работе, является новый способ решения задач выпуклой оптимизации вида\footnote{Важно заметить, что выпуклость задачи понимается в данной статье как выпуклость по совокупности всех переменных, а не просто как выпуклость отдельно по блокам. Последнего может быть недостаточно для того, чтобы исходная задача была выпуклой, и к ней можно было применять обсуждаемые численные методы. Хорошо известными примерами таких задачи (выпуклых по блокам, но не выпуклых в целом) являются задачи, возникающие в рекомендательных системах, например, matrix factorization.} (для наглядности изложение и точности ссылок на литературу рассмотрим задачу без ограничений и композитных членов)
\begin{equation}\label{minmin}
\min_x\min_y f(x,y).
\end{equation}
Если функция $f$ обладает Липшицевым градиентом $\nabla_x f$ по $x$ с константой $L_x$ в $2$-норме и Липшицевым градиентом $\nabla_y f$ по $y$ с константой $L_y$ в $2$-норме, то можно решить задачу \eqref{minmin} с (ожидаемой) точностью $\varepsilon$ по функции за $O\left(\sqrt{\frac{L_xR^2}{\varepsilon}}\right)$  вычислений $\nabla_x f$ и $O\left(\sqrt{\frac{L_yR^2}{\varepsilon}}\right)$  вычислений $\nabla_y f$, где $R^2$~-- расстояние от точки старта до решения в $2$-норме \cite{Kovalev2022}. В частности, один из способов получения данного результата -- это использование ускоренных блочно-покомпонентных методов линейки \cite{gasnikov2016coordinate,Nesterov2017,Allen-Zhu2016}, в которых при выборе параметра метода равным $1/2$ ожидаемое число итераций будет равняться $O\left(\frac{(L_x^{1/2} + L_x^{1/2})R}{\varepsilon^{1/2}}\right)$, а вероятности выбора блока $x$ или $y$ будут равняться, соответственно, $\frac{L_x^{1/2}}{L_x^{1/2} + L_x^{1/2}}$ и $\frac{L_y^{1/2}}{L_x^{1/2} + L_x^{1/2}}$. Откуда и получается приведенный выше результат о расщеплении оракульных сложностей по блокам. Несложно проверить, что данные оценки уже невозможно дальше улучшить \cite{Kovalev2022}.

Если, скажем, по блоку переменных $x$ у нас есть только Липшицевость с константой $M_x$ в $2$-норме, то используя стандартный прием \cite{Nesterov2015} (см. также \cite{gasnikov2016coordinate} по части адаптации к блочно-покомпонентным методам) погружения негладкой задачи в класс гладких задач за счет выбора $L_x = \frac{M_x^2}{2\delta}$ и условия на $\delta$: $N_x\delta = \varepsilon/2$, где $N_x \simeq \sqrt{\frac{L_xR^2}{\varepsilon}}$ -- число вычислений $\nabla_x f$ , получаем, что $N_x \simeq \frac{M_x^2R^2}{\varepsilon^2}$, что также соответствует нижним оценкам \cite{Nemirovski1979}. Отметим, что в \cite{gasnikov2016coordinate} выписано более грубое условие $(N_x +N_y)\delta = \varepsilon/2$, которое, на самом деле, можно переписать так, как это было сделано в данной работе выше.

Таким образом, если в задаче \eqref{minmin}   выпуклая целевая функция $f(x,y)$ негладкая по $x$, но гладкая по $y$  с оракулом $\nabla_x f$, который сильно дешевле, чем оракул $\nabla_y f$, то описанное расщепление оракульных сложностей:  $\sim \varepsilon^{-2}$ вычислений $\nabla_x f$ и $\sim \varepsilon^{-1/2}$ вычислений $\nabla_y f$, дает заметное ускорение по сравнению со стандартным подходом, классифицирующим задачу как негладкую и требующим  $\sim \varepsilon^{-2}$ вычислений $\nabla_x f$ и $\nabla_y f$. 

К сожалению, на рассмотренной в данной статье задаче поиска равновесия в двухстадийной транспортной модели описанный эффект расщепления не проявляется в полной мере, потому что не верно предположение: ``оракул $\nabla_x f$, сильно дешевле, чем оракул $\nabla_y f$''. В нашем случае получается ровно наоборот. Однако, ситуация может поменяться, если для вычисления $\nabla_x f$ можно эффективно (эффективнее, чем для вычисления $\nabla_y f$) использовать GPU/CPU параллелизацию.

Описанный выше подход с изначально заданными $L_x$ и $L_y$ можно строго теоретически обосновать. Проблемы начинаются при обосновании адаптивных (универсальных) версий данного подхода, в которых по блоку $x$ метод должен сам настраиваться на текущий уровень гладкости целевой функции по этому блоку переменных. Возникновение проблемы обусловлено зависимостью вероятностей выбора блоков от адаптивно меняющихся параметров $L_x$ и $L_y$. При этом интересно заметить, что проблема связана именно с выбранным классом ускоренных блочно-покомпонентных методов, обеспечивающих расщепление оракульных сложностей. Для других покомпонентных методов адаптивная настройка на константы Липшица градиента по блокам вполне возможна \cite{Nesterov2012}. В данной работе, рассматривались различные практические способы решения отмеченной теоретической проблемы. В перспективе интересно было бы попробовать получить здесь и теоретическое обосновании полноценной возможности использования конструкции универсальных методов. 

Дополнительно, для демонстрации эффективности покомпонентного метода на описанном классе задач была введена тестовая задача минимизации функции
\newcommand{\xy}{\begin{pmatrix}
       x\\y 
    \end{pmatrix}}
\begin{align*}\label{eq:test_func}
    f(x, y) &= \gamma\ln \left(\sum_{k=1}^m\exp\left({\frac{A_k^T x - b_k}{\gamma}}\right)\right) 
    \\&+
    \frac{1}{2}(x~y) B\xy,\numberthis
\end{align*}
где $B \in \R^{210}$~-- симметричная положительно полуопределённая матрица, $A_k$~-- $k$-ый столбец матрицы $A \in \R^{10 \times 100}$.

Выбором размерность переменных $x \in \R^{10}$,  $y \in \R^{200}$ достигалась разница сложности вычисления градиентов по разным блокам переменных.

Задача решалась при достаточно маленьком $\gamma = 10^{-3}$ и небольшом максимальном собственном числе матрицы $B$, $\lambda_{\max}(B) = 10^{-1}$, а все сингулярные числа матрицы $A$ были равны нулю.
В таком случае, по переменной $x$ задача менее гладкая, то есть имеет большую константу Липшица градиента ($L_x = \lambda_{\max}(B) + 1/\gamma,~L_y = \lambda_{\max}(B)$).
В то же время, из-за существенного различия в размерности переменных $x$ и $y$, вычисление градиента по $x$ требует намного меньше арифметических операций.
В этом случае применение ACRCD* позволяет <<расщепить>> оракульные сложности, делая меньше вычислений градиента по <<дорогому>>, но гладкому блоку переменных $y$ при решении в целом <<негладкой>> задачи (по совокупности переменных).

Для сравнения запускались описанные выше алгоритмы ACRCD* \cite{ACRCD} и USTM \cite{GasinkovNesterov2016}. В качестве метрики, оценкивающей качество приближенного решения, использовалась обратная норма градиента:
\begin{align}\label{eq:metric_test}
    \left\|  \nabla f(x, y) \right\| ^ {-1}.
\end{align}

Графики сходимости изображены на Рисунке ~\ref{fig:test_convergence}.
Для ACRCD* приведены значения метрики в зависимости от количества обращений к каждому из оракулов для визуализации расщепления оракульной сложности. 
Сходимость по норме градиента существенно немонотонная, поэтому для удобства визуализации траетории сходимости ACRCD* были сглажены: точки, изображенные на графике, соответствуют  значениям метрики, усреднённым по соседним 30 итерациям.

Результаты экспериментов показали что метод эффективен на данной задаче -- достичь одного и того же значения метрики можно за значительно меньшее количество обращений к <<дорогому>> оракулу. Это означает что с помощью данного метода для описанных задач достигается общее ускорение по времени работы алгоритма.

\begin{figure*}
    \centering
        \begin{tikzpicture}
    
    \begin{axis}
    [
        cycle list name=exotic,
        xlabel={Метрика \eqref{eq:metric_test}},
        ylabel={Оракульные вызовы},
        xmin=0, xmax=200000,
        width=5cm,
        height=12cm,
        ymin=0, ymax=20000,
        xmode=log,
        ymode=log,
        width=400,
        legend pos=south east,
        ymajorgrids=true,
        grid style=dashed,
    ]

    \addplot table [x=x, y=y, col sep=comma] 
    {test_problem_first_oracle.csv};
    \addplot table [x=x, y=y, col sep=comma] {test_problem_second_oracle.csv};
    \addplot table [x=x, y=y, col sep=comma] 
    {test_problem_ustm.csv};
    
    \legend{
    {ACRCD* - вызовы $\nabla_x f(x, y)$},
    {ACRCD* - вызовы $\nabla_y f(x, y)$},
    {USTM - вызовы $\nabla f(x, y)$},
    }
    \end{axis}
    
    \end{tikzpicture}
    \caption{Сравнение алгоритмов на задаче минимизации функции \eqref{eq:test_func}}
    \label{fig:test_convergence}
\end{figure*}

\section{Заключение}
В данной работе рассматривается задача поиска равновесия в двухстадийной модели транспортных потоков. Задача поиска равновесия сводится к задаче выпуклой минимизации со структурой $\min\min$, где целевая функция обладает существенно разными свойствами по разным группам переменных. Теоретическое исследование показало, что наиболее эффективные методы решения задач такого типа со структурой (в том числе новые, предложенные в данной статье) могут существенно улучшать эффективность стандартных (универсальных) ускоренных процедур. Однако применительно к задаче поиска равновесий в двухстадийной модели, из-за того, что оракул по негладкому блоку переменных сильно дороже оракула по гладкому блоку, получается, что эффект от расщепления оракульных сложностей проявляется не так ярко, как потенциально мог бы проявляться в целом на рассматриваемом классе $\min\min$ задач. Поэтому ключевую роль играют именно численные эксперименты, которые показали, что алгоритм ACRCD*, хоть и является <<физичным>> в контексте транспортной задачи, однако работает хуже, чем обычный USTM и USTM-Sinkhorn!
Принципиальная возможность уменьшения времени вычислений при решении гладко-негладких задач продемонстрирована на синтетическом примере.

Исследование в части теории в разделах 3, 5 было проведено А.\,В. Гасниковым за счет гранта Российского научного фонда (проект № 23-11-00229, \url{https://rscf.ru/project/23-11-00229/}).

Работа Д.\,В. Ярмошика в разделе 2
выполнена  при финансовой поддержке гранта поддержки ведущих научных школ НШ775.2022.1.1.

Практическая часть исследований была выполнена при поддержке ежегодного дохода ФЦК МФТИ ({целевого капитала № 5} на развитие направлений искусственного интеллекта и машинного обучения в МФТИ, \url{https://fund.mipt.ru/capitals/ck5/}).


\newpage

\begin{thebibliography}{99}
\bibitem{baimurzina2019jvm}
	\textit{Баймурзина Д.\,Р. и др.} Универсальный метод поиска равновесий и стохастических равновесий в транспортных сетях~// Журнал вычислительной математики и математической физики. – 2019. – Т. 59. – №. 1. – С. 21-36.

\bibitem
{Wilson1978}
	\textit{Вильсон~А.\,Дж.} Энтропийные методы моделирования сложных систем~// М.: Наука,~--- 1978. 
\bibitem
{gasnikov2013book}
	\textit{Гасников~А.\,В., Кленов~С.\,Л., Нурминский~Е.\,А., Холодов~Я.\,А., Шамрай~Н.\,Б.} Введение в математическое моделирование транспортных потоков. Под ред. А.В. Гасникова с приложениями М.Л. Бланка, К.В. Воронцова и Ю.В. Чеховича, Е.В. Гасниковой, А.А. Замятина и В.А. Малышева, А.В. Колесникова, Ю.Е. Нестерова и С.В. Шпирко, А.М. Райгородского, с предисловием руководителя департамента транспорта г. Москвы М.С. Ликсутова.~// М.: МЦНМО, -- 2013. -- 427 стр., 2-е изд. 
\bibitem
{gasnikov2014matmod}
	\textit{Гасников А.\,В. и др.} О трехстадийной версии модели стационарной динамики транспортных потоков~// Математическое моделирование. – 2014. – Т. 26. – №. 6. – С. 34-70.

\bibitem
{gasnikov2016} 
    \textit{Гасников~А.\,В. и др.} Эволюционные выводы энтропийной модели расчета матрицы корреспонденций // Математическое моделирование~--- 2016.~--- Т. 28, № 4~--- С. 111-124. 
\bibitem
{gasnikov2020posobie}
	\textit{Гасников~А.\,В., Гасникова ~Е.\,В.} Модели равновесного распределения потоков в больших сетях~// М.: УРСС, ~--- 2023.  
\bibitem
{gasnikov2016coordinate}
\textit{Гасников А. В., Двуреченский П. Е., Усманова И. Н.} О нетривиальности быстрых (ускоренных) рандомизированных методов // Труды Московского физико-технического института. – 2016. – Т. 8. – №. 2 (30). – С. 67-100.

\bibitem
{gasnikov2018jvm}
	\textit{Гасников А.\,В., Нестеров Ю.\,Е.} Универсальный метод для задач стохастической композитной оптимизации~// 
	Журнал вычислительной математики и математической физики. – 2018. – Т. 58. – №. 1. – С. 51-68. \\

\bibitem
{kotlyarova2021} \textit{Котлярова Е. В. и др.} Поиск равновесий в двухстадийных моделях распределения транспортных потоков по сети // Компьютерные исследования и моделирование. – 2021. – Т. 13. – №. 2. – С. 365-379.

\bibitem
{kotlyarova2022} 
\textit{Котлярова Е. В.  и др.} Обоснование связи модели Бэкмана с вырождающимися функциями затрат с моделью стабильной динамики // Компьютерные исследования и моделирование. 2022. Т. 14:2. С. 335–342.	

\bibitem
{Nemirovski1979}
\textit{Немировский А. С., Юдин Д. Б.} Сложность задач и эффективность методов оптимизации // М.: Наука, 1979.

\bibitem {Allen-Zhu2016}
 \textit{Allen-Zhu Z. et al.} Even faster accelerated coordinate descent using non-uniform sampling // International Conference on Machine Learning. – PMLR, 2016. – V. 1110--1119.

 \bibitem
 {Boyles2020}
 \textit{Boyles S. D., Lownes N. E., Unnikrishnan A.} Transportation network analysis. Volume I: Static and Dynamic Traffic Assignment, 2020.
 
\bibitem
{DeCea2005}
 \textit{De Cea J., Fernandez J. E., Dekock V. and Soto A.} Solving network equilibrium problems on
multimodal urban transportation networks with multiple user classes // Transport Reviews. 2005. V. 25(3). P. 293--
317.

 \bibitem
 {Evans1976}
 \textit{Evans S. P.} Derivation and analysis of some models for combining trip distribution and
assignment // Transportation Research, 1976. V. 10(1). P. 37-–57.	

\bibitem
{Gasnikova2023_}
\textit{Gasnikova E. et al.} Sufficient conditions for multi-stages traffic assignment model to be the convex optimization problem // arXiv:2305.09069.

\bibitem
{Yarmoshik+Meruza USTM}
\textit{Meruza Kubentayeva et al.} Primal-Dual Gradient Methods for Searching Network
Equilibria in Combined Models with Nested Choice Structure
and Capacity Constraints // arXiv:2307.00427.

\bibitem
{Gasnikova2023}
\textit{Gasnikova E. et al.} An evolutionary view on equilibrium models of transport flows // Mathematics. – 2023. – V. 11. – no. 4. – P. 858.

\bibitem{Kubentaeva2023}
\textit{Kubentaeva M. et al.} Primal-Dual Gradient Methods for Searching Network Equilibria in Combined Models with Nested Choice Structure and Capacity Constraints // arXive:2307.00427

\bibitem{Kovalev2022}
\textit{Kovalev D., Gasnikov A., Malinovsky G.} An Optimal Algorithm for Strongly Convex Min-min Optimization // arXiv:2212.14439. 


\bibitem{Guminov2019}
\textit{Guminov S. et al.} Accelerated alternating minimization~// ICML 2021

  \bibitem{Nesterov2012}
 \textit{Nesterov Y.} Efficiency of coordinate descent methods on huge-scale optimization problems // SIAM Journal on Optimization. – 2012. – V. 22. – no. 2. – P. 341--362.

 \bibitem{Nesterov2005}
	\textit{Nesterov Y.} Smoothing technique and its applications in semidefinite optimization //Mathematical Programming. – 2007. – V. 110. – no. 2. – P. 245--259.
 

 \bibitem{Nesterov2015}
	\textit{Nesterov Y.} Universal gradient methods for convex optimization problems~// Mathematical Programming. – 2015. – V. 152. – no. 1-2. – P. 381--404.

 \bibitem{GasinkovNesterov2016}
    \textit{Gasnikov--Nesterov} A universal method for stochastic problems
    composite optimization.~// arxiv preprint arXiv:1604.05275. - 2016




\bibitem
{Nesterov2003}
	\textit{Nesterov Y., De Palma A.} Stationary dynamic solutions in congested transportation networks: summary and perspectives~// Networks and spatial economics. – 2003. – Т. 3. – No. 3. – P. 371-395.

\bibitem
{Nesterov2017}
 \textit{Nesterov Y., Stich S.} Efficiency of the accelerated coordinate descent method on structured optimization problems // SIAM Journal on Optimization. – 2017. – V. 27. – no. 1. – P. 110--123.

\bibitem
{Ortuzar2002}
	\textit{Ortúzar J.\,D., Willumsen L.\,G.} Modelling transport. John Wiley and Sons~// West Sussex, England. – 2002.
	
\bibitem
{Patriksson2015}
	\textit{Patriksson M.} The traffic assignment problem: models and methods. – Courier Dover Publications, 2015.
	
\bibitem
{Peyre2019}
	\textit{Peyré G.,  Cuturi M.} Computational Optimal Transport: With Applications to Data Science~// Foundations and Trends® in Machine Learning. – 2019. – V. 11. – №. 5-6. – P. 355-607.

\bibitem{ACRCD}
    \textit{А.В. Гасников, П.Е. Двуреченский, Усманова И.Н.} О нетривиальности быстрых (ускоренных) рандомизированных методов.~// arxiv preprint arXiv:1508.02182 - 2015



	
\bibitem
{Stabler2020}
	\textit{Stabler B., Bar-Gera H., Sall E.} Transportation Networks for Research Core Team. Transportation Networks for Research. Accessed Month, Day, Year. [Electronic resource]: 
	\url{https://github.com/bstabler/TransportationNetworks} (accessed 16.02.2021).

\bibitem{Tupitsa2022}
 \textit{Tupitsa N. et al.} Numerical Methods for Large-Scale Optimal Transport // arXiv:2210.11368.


\bibitem{GasnikovUGD}
 \textit{Гасников, А. В.} Современные численные методы оптимизации. Метод универсального градиентного спуска //  arXiv:1711.00394.
 
\bibitem{Aliev2015}
\textit{Алиев, А.С., Мазурин, Д.С., Максимова, Д.А. и Швецов, В.И., 2015}. Структура комплексной модели транспортной системы г. Москвы. // Труды Института системного анализа Российской академии наук, 65(1), pp.3-15.

\bibitem{Code}
    \textit{Ссылка на репозиторий с исходным кодом вычислительных экспериментов}
    \url{https://github.com/Lareton/transport_network_optimization}

	
\end{thebibliography}
\end{document}